\documentclass[graybox]{svmult}

\usepackage{mathptmx}       

\usepackage{helvet}         
\usepackage{courier}        

\usepackage{graphicx}              

\usepackage{siunitx}
\usepackage{amsbsy}
\usepackage{amsmath,bm}
\usepackage{amssymb}
\usepackage{amsfonts}

\usepackage{multirow}

\usepackage[bottom]{footmisc}
\usepackage{cite}

\usepackage{url}

\def\IR{{\mathbb R}}

\newcommand{\bA}{{\textbf A}}
\newcommand{\bB}{{\textbf B}}
\newcommand{\bC}{{\textbf C}}
\newcommand{\bD}{{\textbf D}}

\newcommand{\bY}{{\textbf Y}}

\newcommand{\bN}{{\textbf N}}
\newcommand{\bI}{{\textbf I}}

\newcommand{\bO}{{\textbf O}}

\newcommand{\bW}{{\textbf W}}
\newcommand{\bR}{{\textbf R}}

\newcommand{\bx}{{\textbf x}}
\newcommand{\by}{{\textbf y}}
\newcommand{\bu}{{\textbf u}}
\newcommand{\bV}{{\textbf V}}
\newcommand{\bU}{{\textbf U}}

\newcommand{\bq}{{\textbf q}}

\newcommand{\bw}{{\textbf w}}

\newcommand{\cO}{ {\cal O} }
\newcommand{\cH}{ {\cal H} }

\newcommand{\cS}{ {\cal S} }
\newcommand{\cW}{ {\cal W} }

\newcommand{\bSigma}{\boldsymbol{\Sigma}}

\newcommand{\cR}{ {\cal R} }

\begin{document}

\title*{Bilinear realization from i/o data with NN}
\author{D.~S.~Karachalios
\and
I.~V.~Gosea
\and
K.~Kour
\and
A.~C.~Antoulas
}
\institute{D.~S. Karachalios, I.~V. Gosea, K. Kour \at Max Planck Institute for Dynamics of Complex Technical Systems, Magdeburg, Germany, \email{{karachalios, gosea, kour}@mpi-magdeburg.mpg.de}
\and A.~C. Antoulas \at Max Planck Institute for Dynamics of Complex Technical Systems, Magdeburg, Germany,\at Rice University, Electrical and Computer Engineering Department, Houston, TX 77005, USA,\at Baylor College of Medicine, Houston, TX 77030, USA, \email{aca@rice.edu}}

%
%
\maketitle

\abstract*{We present a method that connects a well-established nonlinear (bilinear) identification method from time-domain data with neural network (NNs) advantages. The main challenge for fitting bilinear systems is the accurate recovery of the corresponding Markov parameters from the input and output measurements. Afterward, a realization algorithm similar to that proposed by Isidori can be employed. The novel step is that NNs are used here as a surrogate data simulator to construct input-output (i/o) data sequences. Then, classical realization theory is used to build a bilinear interpretable model that can further optimize engineering processes via robust simulations and control design.}
\vspace{-40mm}
\abstract{We present a method that connects a well-established nonlinear (bilinear) identification method from time-domain data with neural network (NNs) advantages. The main challenge for fitting bilinear systems is the accurate recovery of the corresponding Markov parameters from the input and output measurements. Afterward, a realization algorithm similar to that proposed by Isidori can be employed. The novel step is that NNs are used here as a surrogate data simulator to construct input-output (i/o) data sequences. Then, classical realization theory is used to build a bilinear interpretable model that can further optimize engineering processes via robust simulations and control design.}
\vspace{-5mm}
\section{Introduction}\label{sec:intro}
\vspace{-5mm}
Evolutionary phenomena can be described formally as continuous dynamical models with partial differential equations (PDEs). The continuous nature of these physical models is equipped with analytical results for efficient discrete approximation in space and time. In particular, methods such as finite elements or differences bridge the continuous analytical laws of the physical world with computational science that has discrete nature \cite{ACA05}. On the other hand, data science allows model discovery when the identification feature is under consideration \cite{ABG20}. The quantification of the above equivalences, in combination with the stochastic nature governing real-world applications, aims to explain the digital twin. Spatial discretization of PDEs in many cases results in a continuous w.r.t. time nonlinear system of ordinary differential equations (ODEs) described by the operators $(\mathcal{F},~\mathcal{G})$ and it can be approximated with Carleman linearization (e.g., bilinear system form) \cite{Car32}.
\begin{equation}\label{eq:nonlinBil}
   \bSigma:\left\{\begin{aligned}
    \dot{\bx}(t)&=\mathcal{F}(\bx(t),u(t))\\
    y(t)&=\mathcal{G}(\bx(t),u(t))
    \end{aligned}\right.\xrightarrow[\bSigma\approx\bSigma_{bil}]{\text{Carleman}}\left\{\begin{aligned}
    \dot{\bx}(t)&=\mathcal{A}\bx(t)+\mathcal{N}\bx(t)u(t)+\mathcal{B} u(t)\\
    y(t)&=\mathcal{C}\bx(t)+\mathcal{D} u(t),~\bx_0=\textbf{0},~t\geq 0.
    \end{aligned}\right.
\end{equation}
If the original system has dimension $N$, since Carleman linearization \cite{Car32} preserves up to the quadratic term $\bx(t)\otimes\bx(t)$\footnote{$\otimes$: Kronecker product}, the dimension of the resulting bilinear system $(\mathcal{A},~\mathcal{N},~\mathcal{B},~\mathcal{C},~\mathcal{D})$ increases to $n=N^2+N$.

Data-driven methods can be classified into two general classes. The first provides prediction via regression techniques such as neural networks (NN) from machine learning (ML), while the second has its roots in system theory and allows model discovery \cite{ABG20,AGI16}. Generally, the NNs are sensitive to parameter tuning and lack the model interpretability because of the inherent "black-box" structure \cite{Schilders2008}, while the latter constructs interpretable models and are capable of explaining the hidden dynamics. ML models learn the features through a composition of a couple of nonlinear functions from data during the training process. Hence, by using data points for the training, the prediction would be expressed as a function of these data points. Until recently, the ML and system identification (SI) techniques were developed independently. But in recent years, there has been a great effort to establish a common ground \cite{LJUNG2011}.

The authors in \cite{Favoreel1999} have generalized the linear subspace identification theory to bilinear systems. In chemical processes, the controls are flow rates and from first principles (mass and heat balances) these will appear in the system equations as products with the state variables, therefore the bilinear equation has the physical form $\footnotesize M\dot{\bx}=\sum_i\bq_i\bx_i-\sum_m \bq_m\bx_m,~\bq\text{(inputs)},~\bx\text{(state concentration)}$. The authors in \cite{Ramos2009} were able to construct bilinear systems with white noise inputs based on an iterative deterministic-stochastic subspace approach. The author in \cite{Juang2005} makes use of the linear model properties of the bilinear system when subjected to constant input. Constant inputs can transform the bilinear model to an equivalent linear model \cite{GKApamm2021}.

In Sec.\;(\ref{sec:ioBil}), we introduce the bilinear realization theory by explaining in detail the data acquisition procedure to compute the bilinear Markov parameters that will enter the bilinear Hankel matrix. Further, we present a concise algorithm that can achieve bilinear identification, detailed by two examples. In Sec.\;(\ref{sec:NN}), we train a neural network with a single i/o data sequence to mimic the unknown simulator, and we combine it with the bilinear realization theory. As a result, we could construct a bilinear model from a single i/o data with better fit-performance compared with another state-of-the-art bilinear SI approach.       
\vspace{-5mm}
\section{The bilinear realization framework}\label{sec:ioBil}
\vspace{-5mm}
In the case of linear systems, Ho and Kalman \cite{HOKALMAN} have provided the mathematical foundations for achieving linear system realization from i/o data. In the nonlinear case and towards the exact scope of identifying nonlinear systems, Isidori in \cite{Isidori1973} has extended these results for the bilinear case, and Al Baiyat in \cite{AlBaiyat2004} has provided an SVD-based algorithm.  

Time discretization as in \cite{morBenBD11} of the single-input single-output (SISO) bilinear system (\ref{eq:nonlinBil}) of state dimension $n$ with sampling time $\Delta t$, results to fully discrete models defined at time instances given by $0<\Delta t<2\Delta t<\cdots<k\Delta t$, with $\bx_c(k\Delta t)=\bx_k$ and $u(k\Delta t)=u_k$ for $k=0,\ldots,m-1$
\begin{equation}\label{eq:discBilSysZ}\footnotesize
    \bSigma_{\text{disc}} :\left\{\begin{aligned}
\bx_{k+1}&=\bA\bx_k+\bN\bx_k u_k+\bB u_k,\\
y_k&=\bC\bx_k + \bD u_k,~\bx_0=\textbf{0}.
\end{aligned}\right.
\end{equation}
The discrete-time system in (\ref{eq:discBilSysZ}) has state dimension $n$, so, $\bx\in\IR^n$ and the operators have dimensions $\bA,\bN\in\IR^{n\times n},~\bB,\bC^T\in\IR^{n}$, and $\bD\in\IR$. We assume homogeneous initial conditions and a zero $\bD$ term. The forward Euler is the only numerical scheme that preserves the bilinear structure in the discrete set-up with the cost of conditional stability. Moreover, a more sophisticated scheme can interpolate exactly the continuous model at the sampling points in \cite{Dunoyer2010} but is restricted to only a subclass of bilinear systems. Therefore, one good choice in terms of stability is the backward Euler scheme from \cite{BBD2011}, which preserves the bilinear structure asymptotically.
\begin{equation}\label{eq:discCont}
    \phi:~\bA=(\bI-\Delta t\mathcal{A})^{-1},~\bN=\Delta t(\bI-\Delta t\mathcal{A})^{-1}\mathcal{N},~\bB=\Delta t(\bI-\Delta t\mathcal{A})^{-1}\mathcal{B},~\bC=\mathcal{C}.
\end{equation}
Note that the transformation $\phi$ is invertible:  $\bSigma_{b}^{c}:(\mathcal{A},\mathcal{N},\mathcal{B},\mathcal{C})\leftrightarrow\bSigma_{b}^{d}:(\bA,\bN,\bB,\bC)$.
\begin{definition}\label{def:reach} The reachability matrix $\cR_{n}=\left[\begin{array}{ccc}
    \bR_{1} & \cdots & \bR_{n}  \\
\end{array}\right]$ is defined recursively from the following relation: $\bR_j=\left[\begin{array}{cc}
  \bA\bR_{j-1} & \bN\bR_{j-1}\end{array}\right],~j=2,\ldots n,~\bR_1=\bB$.
\end{definition}
Then, the state space of the bilinear system is spanned by the states reachable from the origin if and only if $\text{rank}(\cR_n)=n$.
\begin{definition}\label{def:observ} The observability matrix $\cO_{n}=\left[\begin{array}{ccc}
    \bO_{1} & \cdots & \bO_{n}  \\
\end{array}\right]^T$ is defined recursively from the following relation: $\bO_j^T=\left[\begin{array}{cc}
  \bO_{j-1}\bA & \bO_{j-1}\bN\end{array}\right]^T,~j=2,\ldots n,~\bO_1=\bC$.
\end{definition}
Then the state space of the bilinear system is observable iff $\text{rank}(\cO_n)=n$. The following Def.\;(\ref{def:input}) will allow a concise representation of the i/o relation.
\begin{definition}\label{def:input} 
$\bu_j(h)=\left[\begin{array}{c}
    \bu_{j-1}(h) \\ \bu_{j-1}(h)u(h+j-1) 
\end{array}\right],~j=2,\ldots,~\bu_1(h)=u(h)$.
\end{definition}
Let $\{\bw_1,\bw_2,\ldots,\bw_j,\ldots\}$ be an infinite sequence of row vectors, in which $\bw_j\in\IR^{1\times 2^{j-1}}$ and is defined recursively as follows $\bw_j=\bC\bR_{j},~j=1,2,\ldots$;
\vspace{-5mm}
\subsection{The bilinear Markov parameters}
\vspace{-5mm}
Bilinear Markov parameters are encoded in the $\{\bw_j\}$-vectors for $j=1,\ldots,n$. These are invariant quantities of the bilinear system in connection with the input-output relation. After making use of Def.\;(\ref{def:input}), we can write 
\begin{equation}\label{eq:LSunder}\footnotesize
   \underbrace{\left[\begin{array}{c}
        y_1  \\
        y_2 \\
        \vdots\\
        y_k
   \end{array}\right]}_{\bY}=\underbrace{\left[\begin{array}{cccc}
        \bu_{1}^T(0) & 0 & \cdots & 0  \\
        \bu_{1}^T(1) & \bu_{2}^T(0) & \cdots & 0  \\
        \vdots & \vdots & \ddots & \vdots  \\
        \bu_{1}^T(k-1) & \bu_{2}^T(k-2) & \cdots & \bu_{k}^T(0)
    \end{array}\right]}_{\bU}\cdot\underbrace{\left[\begin{array}{c}
    \bw_1^T \\
    \bw_2^T \\
    \vdots \\
    \bw_k^T\end{array}\right]}_{\bW},
\end{equation}
where the dimensions are: $\bY\in\IR^{k\times1}$, $\bU\in\IR^{k\times m}$, and $\bW\in\IR^{m\times 1}$. 
The least squares problem filled out with $k$ time steps will remain under-determined $\forall k\in\{2,3,\ldots\}$ as long as $m=2^k-1$ bilinear Markov parameters will be activated. Thus, we must deal with $k$ equations and $2^k-1$ unknowns. Solving an under-determined system is not impossible, but the solutions are infinite, and regularization schemes cannot easily lead to identification. Therefore, one way of uniquely identifying the bilinear Markov parameters and determining the $\bW$ solution vector can be achieved by solving a coupled least squares system after applying several simulations to the original system. 

To determine $(2^k-1)$ parameters, the column rank of the matrix $\bU$ should be full. That can be accomplished by augmenting rows underneath to the matrix $\bU$ up to the point where the new extended matrix $\hat{\bU}$ will be at least square as the result of concatenating several $\bU$ matrices that are coming from the different simulations with different inputs. Thus, we need at least $2^{k-1}$ independent simulations from the original system. That is exactly the bottleneck expected for nonlinear identification frameworks that deal with time-domain data. Later on, we will relax this condition in a novel way by using NNs. Equation (\ref{eq:getMarkov}) describes the coupled linear least squares system with $d=2^{k-1}$ independent simulations that can provide the unique solution $\cW$ with the bilinear Markov parameters. 
\begin{equation}\label{eq:getMarkov}\footnotesize
    \underbrace{\left[\begin{array}{ccc}
         \bY_{1}&\cdots&\bY_{d}
    \end{array}\right]^T}_{\hat{\bY}}=\underbrace{\left[\begin{array}{ccc}
          \bU_{1}&\cdots\bU_{d} 
    \end{array}\right]^T}_{\hat{\bU}}\cdot\cW
\end{equation}
Hence, we repeat the simulation $d$ times, and each time we get $k$ equations, with the $i^{\text{th}}$ simulation to be $\bY_i=\left[\begin{array}{cccc} y_1^{(i)} & y_2^{(i)} & \cdots & y_k^{(i)} \end{array}\right]^T$ and accordingly for the $\bU_i$, the real matrix $\hat{\bU}$ has dimension $2^k\times(2^k-1)$. The $\hat{\bU}$ matrix results after concatenating all the lower triangular matrices with full column rank. To enforce that $\hat{\bU}$ will also have full column rank, one choice is to use white inputs (sampled from a Gaussian distribution) for the simulations. Using of white inputs is widespread for SI. Still, in that case, a careful choice of deterministic inputs can make the inversion exact and capable of recovering the bilinear Markov parameters. For the solution holds: $\text{rank}(\hat{\bU})=2^k-1,~\text{so, the unique solution is:}~\cW=\hat{\bU}^{-1}\hat{\bY}\in\IR^{2^k-1}$. The vector $\cW$ contains the $2^k-1$ bilinear Markov parameters. As we have computed the bilinear Markov parameters, a generalized Hankel matrix can be constructed.
\vspace{-8mm}
\subsection{The bilinear Hankel matrix}
\vspace{-5mm}
The bilinear Hankel matrix is the product between the observability and the reachability matrices. The bilinear Hankel matrix is denoted with $\cH_b$ and is defined as the product of the following two infinite matrices $\cO,~\cR$
\begin{equation}\label{eq:bilHankel}\footnotesize
\cH_b=\cO\cR=\left[\begin{array}{c}
         \bC  \\
         \bC\bA \\
         \bC\bN \\
         \vdots
    \end{array}\right]\left[\begin{array}{cccc}
    \bB & \bA\bB & \bN\bB & \cdots
    \end{array}\right]=\left[\begin{array}{cccc}
        \bC\bB & \bC\bA\bB & \bC\bN\bB & \cdots  \\
        \bC\bA\bB & \bC\bA^2\bB & \bC\bA\bN\bB & \cdots \\
        \bC\bN\bB & \bC\bN\bA\bB & \bC\bN^2\bB & \cdots \\
        \vdots & \vdots & \vdots & \ddots
    \end{array}\right]
\end{equation}
Equation (\ref{eq:bilHankel}) reveals the connection with the bilinear Markov parameters $\cW=\bC\cR$ that appear in the first row of $\cH_b$. In general, the construction of the bilinear Hankel matrix is described in \cite{Isidori1973} with the partial and completed realization theorems along with the partitions\footnote{$\cS^{A}=\{\text{set of $\cH_b$ columns}:\text{from}~2^m~\text{to}~(3\cdot 2^{m-1}-1),~m=1,2,\ldots\}$, \\$~~~\cS^{N}=\{\text{set of $\cH_b$ columns}:\text{from}~3\cdot2^{m-1}~\text{to}~(2^{m+1}-1),~m=1,2,\ldots\}$.} $\cS^{A},~\cS^{N}$ \cite{AlBaiyat2004}.
\vspace{-8mm}
\subsection{Bilinear realization algorithm}\label{sec:algorithm}
\vspace{-5mm}
\textbf{Input:} Input-output time-domain data from  a system $u\rightarrow\boxed{\bSigma?}\rightarrow y$.\\
\textbf{Output:} A minimal bilinear system $(\bA_r,\bN_r,\bB_r,\bC_r)$ of low dimension $r$ that $\bSigma_r\approx\bSigma$.
\begin{enumerate}
    \item Excite the system $\bSigma$ $k$-times with $\bu_m\sim \mathcal{\bN} (\mu,\sigma)$ and collect $\by_m$, where $k=2^{m-1}$.
    \begin{tabular}{cc}
\hline\\[-3mm]
 1st simulation & $[u_{1}(1)\cdots u_{1}(m)]\rightarrow\boxed{\bSigma}\rightarrow[y_{1}(1)\cdots y_{1}(m)]=\bY_1$, and $\bU_1$ as in Def.\;\Ref{def:input}. \\
$\vdots$ & $\vdots$\\
kth simulation & $[u_{k}(1)\cdots u_{k}(m)]\rightarrow\boxed{\bSigma}\rightarrow[y_{k}(1)\cdots y_{k}(m)]=\bY_k$, and $\bU_k$ as in Def.\;\Ref{def:input}.\\[1mm]\hline
\end{tabular}
\item Identify the $(2^m-1)$ bilinear Markov parameters by solving the system in (\ref{eq:getMarkov}).
\item Construct the bilinear Hankel matrix $\cH_b$ and the sub-matrices $\cS^{\bA},~\cS^{\bN}$.
\item Compute $[\bU,\bSigma,\bV]=\text{SVD}(\cH_b)$ and truncate w.r.t the singular values decay ($r\ll n$) - the reduced/identified bilinear model $(\bA_r,\bN_r,\bB_r,\bC_r)$ is constructed 
\begin{eqnarray}\label{eq:quadrupletA}\footnotesize
 \bA_r&=&\bSigma^{-1/2}\bU^T\cS^{\bA}\bV\bSigma^{-1/2}\\
 \label{eq:quadrupletN}
 \bN_r&=&\bSigma^{-1/2}\bU^T\cS^{\bN}\bV\bSigma^{-1/2}\\
 \label{eq:quadrupletB}
 \bB_r&=&\bSigma^{1/2}\bV^{T}~~\rightarrow~~\text{1st column}\\
 \label{eq:quadrupletC}
 \bC_r&=&\bU\bSigma^{1/2}~~\rightarrow~~\text{1st row}
 \end{eqnarray}
\end{enumerate}

\begin{example}(A toy system)
Let the following bilinear system of order $2$ be 
\begin{equation}\footnotesize
\bA=\left[\begin{array}{cc}
0.9 & 0.0\\
0.0 & 0.8
\end{array}\right],~\bN=\left[\begin{array}{cc}
0.1 & 0.2\\
0.3 & 0.4
\end{array}\right],~\bB=\left[\begin{array}{c}
1.0\\
0.0 
\end{array}\right],~\bC=\left[\begin{array}{c}
1.0\\
1.0
\end{array}\right]^T.
\end{equation}
Applying the algorithm in Sec.\;(\ref{sec:algorithm}), by choosing $m=4$, we can recover $2^m-1=15$ bilinear Markov parameters. The solution of the system in Eq.\;\eqref{eq:getMarkov} is:
\begin{equation*}
    \footnotesize\bW=\left[\begin{array}{ccccccccccccccc} 1.0 & 0.9 & 0.4 & 0.81 & 0.33 & 0.36 & 0.22 & 0.729 & 0.273 & 0.297 & 0.183 & 0.324 & 0.18 & 0.198 & 0.118 \end{array}\right]
\end{equation*}
By reshuffling the vector $\cW$, we can form the $\cH_b$ matrix and the shifted versions $\cS^{\bA},~\cS^{\bN}$ as described above. The Hankel matrix ($3$ rows \& $7$ columns are displayed) along with the shifted versions are $\left\{\cH_b,~\cS^A,~\cS^N\right\}:=$
\begin{equation*}\scriptsize\left\{\left[\begin{array}{ccccccc} 1.0 & 0.9 & 0.4 & 0.81 & 0.33 & 0.36 & 0.22\\ 0.9 & 0.81 & 0.33 & 0.729 & 0.273 & 0.297 & 0.183\\ 0.4 & 0.36 & 0.22 & 0.324 & 0.18 & 0.198 & 0.118 \end{array}\right],~\left[\begin{array}{ccc} 0.9 & 0.81 & 0.33\\ 0.81 & 0.729 & 0.273\\ 0.36 & 0.324 & 0.18 \end{array}\right],~\left[\begin{array}{ccc} 0.4 & 0.36 & 0.22\\ 0.33 & 0.297 & 0.183\\ 0.22 & 0.198 & 0.118 \end{array}\right]\right\}.
\end{equation*}
 In Fig.\;(\ref{fig:Exfig1}), the $3$rd normalized singular value has reached machine precision $\sigma_{3}/\sigma_{1}=5.2501e-17$, that is the criterion for choosing the order of the fitted system (which is minimal, in this case) of the underlying bilinear system. Therefore, we construct a bilinear model of order $r=2$ and the realization obtained is equivalent to the original (minimal) one, up to a coordinate (similarity) transformation. Other ways of constructing reduced models from Hankel$\subset$Loewner matrices can be obtained with the CUR (cross approximations based) decomposition scheme as in \cite{KGAhandbook}.
 \begin{equation}
\scriptsize\bA_r=\left[\begin{array}{cc} 0.89394 & 0.11305\\ 0.0050328 & 0.80606 \end{array}\right],~\bN_r=\left[\begin{array}{cc} 0.41116 & -0.2281\\ -0.24782 & 0.088841 \end{array}\right],~\bB_r=\left[\begin{array}{c} -1.0001\\ -0.053577 \end{array}\right],~\bC_r=\left[\begin{array}{c} -1.0001\\ 0.0040101 \end{array}\right]^T.
\end{equation}
\vspace{-10mm}
\begin{figure}[h!]
    \centering
    \includegraphics[scale=0.15]{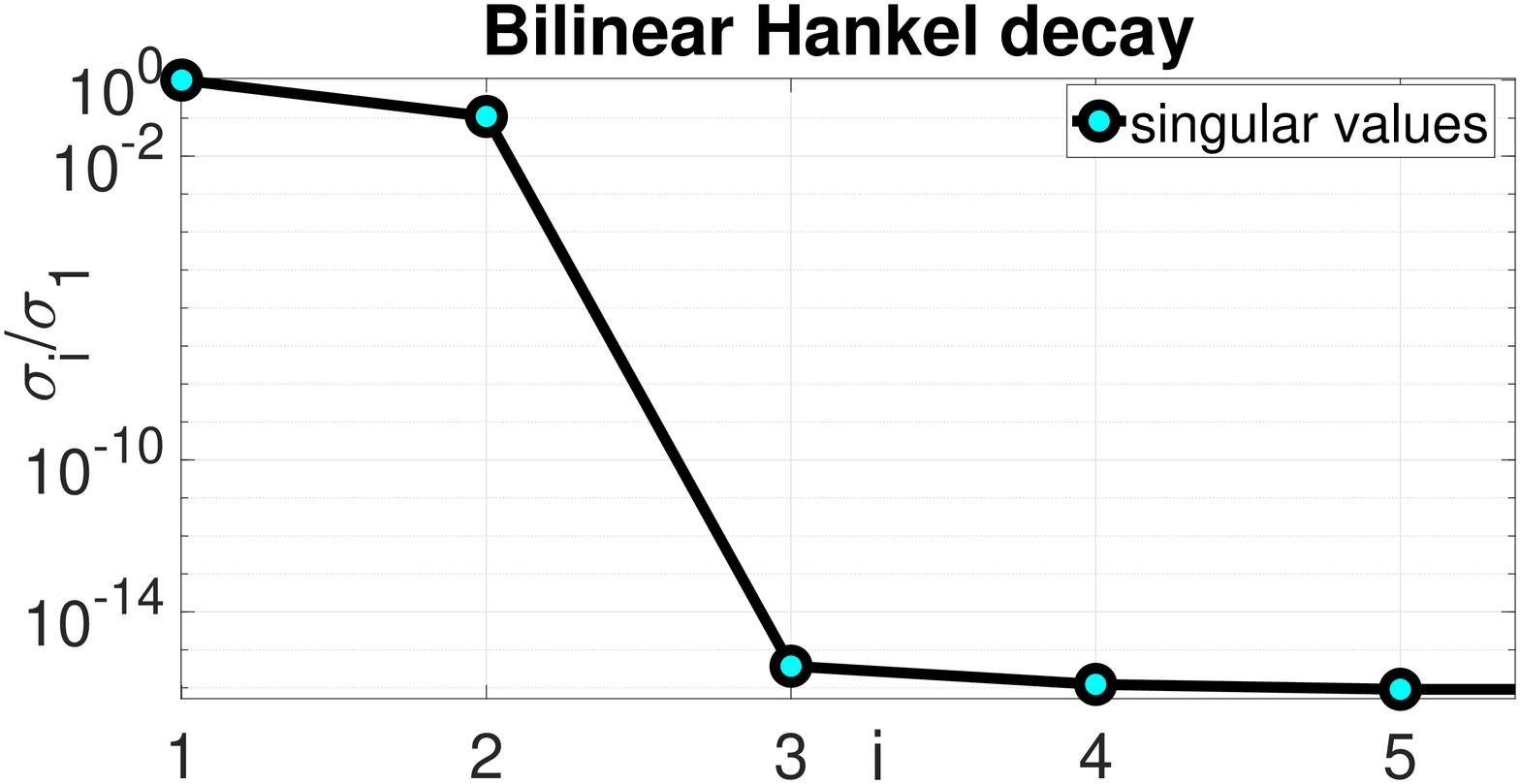}
    \includegraphics[scale=0.15]{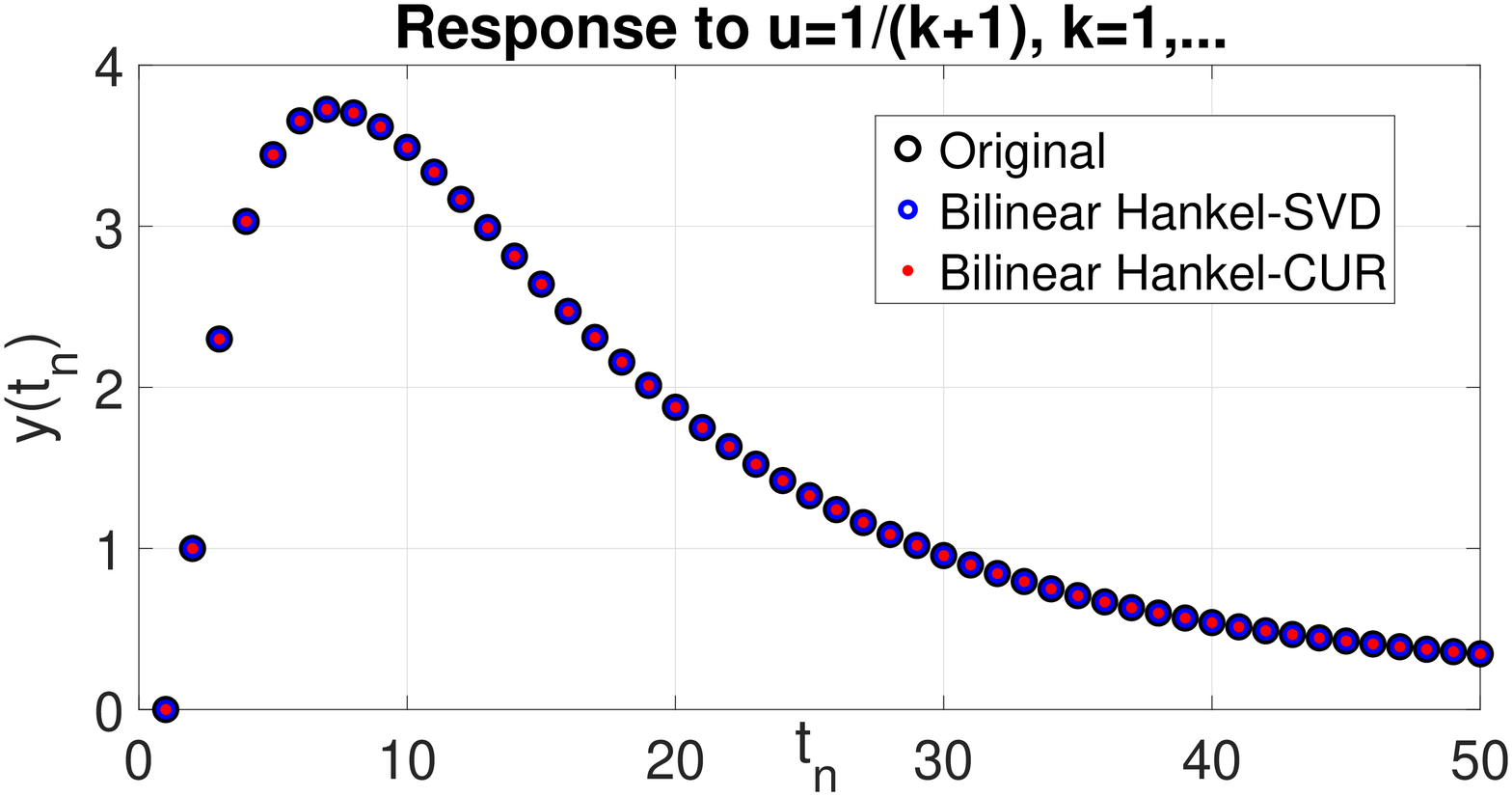}
    \caption{On the left figure, the singular value decay of the bilinear Hankel matrix is depicted. On the right figure, the input response $u_k=1/(k+1),~k=0,1,\ldots$ certifies that all models are equivalent.}
    \label{fig:Exfig1}
\end{figure}
\end{example}
\vspace{-10mm}
\begin{example}(The viscous Burgers' equation example)
Following \cite{BBD2011} after spatial semi-discretization and the Carleman linearization technique, yields a bilinear system of dimension $n=30^2+30=930$. The viscosity parameter is $\nu=0.1$, the sampling time is $\Delta t=0.1$ and with $2^{m-1}=512$ independent random inputs with length $m=10$ each, we construct a database of $5,120$ points. Solving (\ref{eq:getMarkov}), we get the bilinear Markov parameters, and the bilinear Hankel is constructed. On the left pane of Fig.\;(\ref{fig:Burgers}), the bilinear Hankel singular value decay captures the nonlinear nature of the Burgers' equation, whereas, on the other hand, the linear Hankel framework captures only the linear minimal response. It is evident on the right pane of Fig.\;(\ref{fig:Burgers}) after using the inverse transformation $\phi$ from (\ref{eq:discCont}), the reduced continuous-time bilinear model of order $r=18$ performs well, yielding an error $O(10^{-5})$ where at the same time, the linear fit is off.
\vspace{-5mm}
\begin{figure}[h!]
\centering
    \includegraphics[scale=0.17]{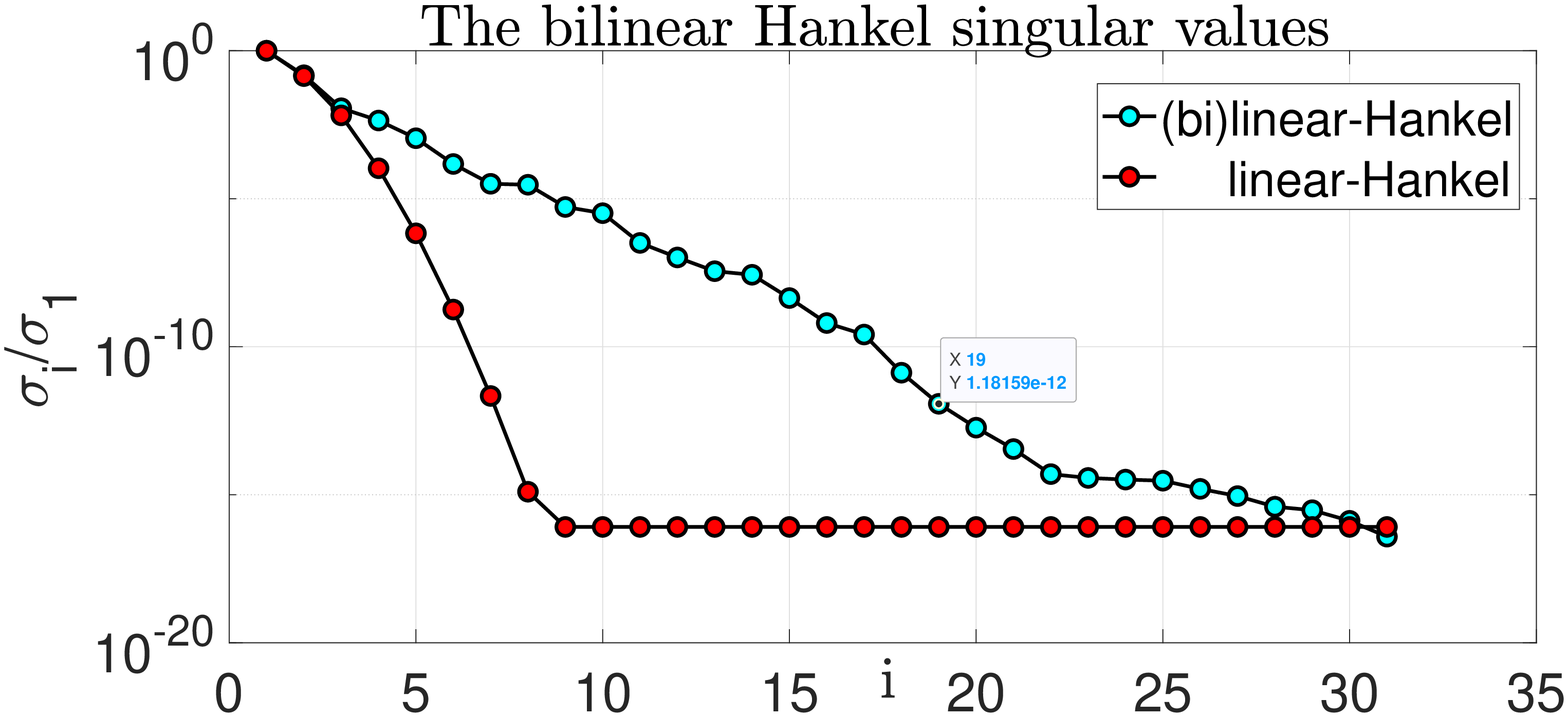}
    \includegraphics[scale=0.17]{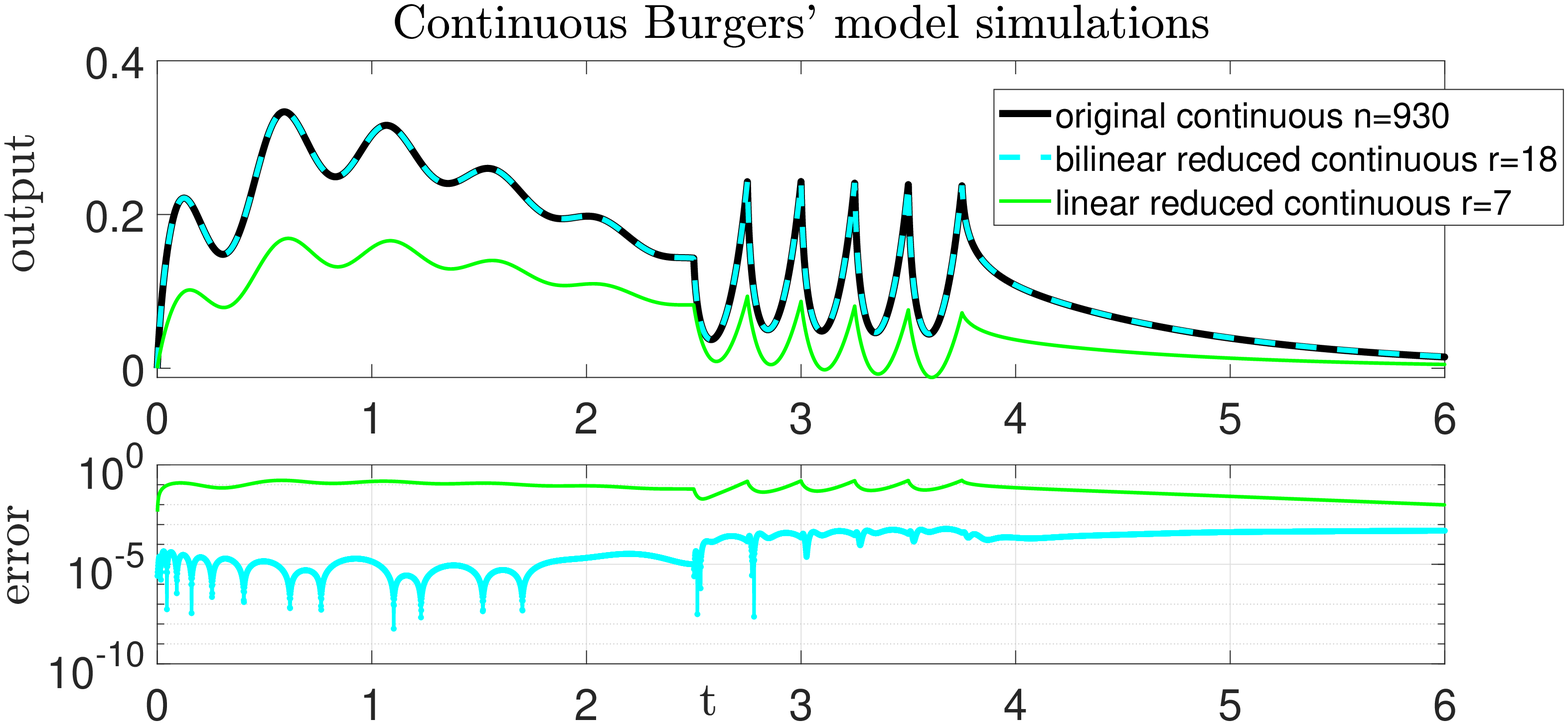}
\caption{\textbf{Left pane}: The recommended reduced bilinear model with $\Delta t=0.1$ is of order $r=18$ where $\sigma_{19}/\sigma_1=1.18\cdot10^{-12}$. \textbf{Right pane}: $u_1=(1+\cos(2\pi t))e^{-t},t\in[0,2.5],u_2=2\text{sawtooth}(8\pi t),t\in[2.5,3.75],u_3=0$, it is compared with a continuous bilinear identification method based on the Loewner framework in both frequency and time domain approaches \cite{Karachalios2021,AGI16}.}
\label{fig:Burgers}
\end{figure}
\vspace{-10mm}
\end{example}
\vspace{-8mm}
\section{From a single data sequence to bilinear realization}\label{sec:NN}
\vspace{-5mm}
To achieve bilinear realization as in \cite{Isidori1973} a repetitive data assimilation simulation in the time domain is required. In many cases, the data from a simulated system are available as a single i/o sequence \cite{Ramos2009}. By using NARX-net based model, in case of a single experiment, in a real engineering environment, the expensive repetitive simulations can be avoided. These models learn from a unique data sequence and are capable of predicting the output behavior under different excitations. That is precisely where NARX-net model architecture will play the role of a surrogate simulator. Then, by constructing a NN-based model \cite{Billings} and combining the realization theory in  \cite{Isidori1973}, a state-space bilinear model can be constructed as in (\ref{eq:discBilSysZ}). The use of a state-space model, which relies on the classical nonlinear realization theory with many known results (especially on bilinear systems and in the study direction of stability, approximation and control) is beneficial in comparison with the NARX. 
\begin{example}(Heat exchanger)
The process is a liquid-saturated steam heat exchanger, where water is heated by pressurized saturated steam through a copper tube. The input variable is the liquid flow rate, and the output variable is the outlet liquid temperature. The sampling time is $1$s, and the number of samples is $4,000$. Further details can be found in \cite{BITTANTI1997135}, and data can be downloaded from the database to identify systems (DaISy):~\url{https://homes.esat.kuleuven.be/~tokka/daisydata.html}.  
\vspace{-6mm}
\begin{figure}[h!]
    \centering
    \includegraphics[scale=0.2]{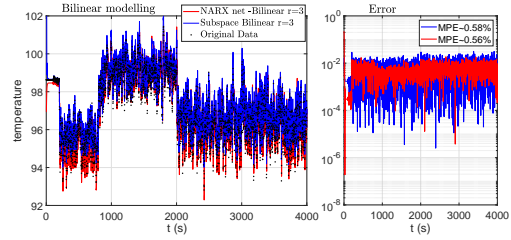}
    \caption{Comparison and model fit of the proposed NARX-net bilinear model (\ref{eq:daisyfit}) with the subspace method from \cite{Ramos2009} for the same reduced order $(r=3)$.}
    \label{fig:DaisyFit}
    \vspace{-8mm}
\end{figure}
\vspace{-2mm}
\begin{equation}\label{eq:daisyfit}\tiny
    \left\{\begin{aligned}
    \dot{\bx}(t)&=\left[\begin{array}{ccc} 0.9164 & 0.09167 & -0.1847\\ -0.2663 & -0.1515 & 0.1232\\ -0.07227 & 0.4778 & 0.3571 \end{array}\right]\bx(t)+\left[\begin{array}{ccc} 0.02717 & 0.5169 & 0.5555\\ -0.09674 & 0.5467 & 0.5696\\ 0.1878 & -0.06846 & -1.981 \end{array}\right]\bx(t)u(t)+\left[\begin{array}{c} 2.9063\\ 2.909\\ -0.16088 \end{array}\right]u(t)+\left[\begin{array}{c} -1.073\\ -1.074\\ 0.05938 \end{array}\right],\\
    y(t)&=\left[\begin{array}{ccc} -0.7852 & 0.7794 & -0.05203 \end{array}\right]\bx(t)+96.9358,~\bx(0)=\textbf{0},~t\geq 0.
    \end{aligned}\right.
\end{equation}
Figure (\ref{fig:DaisyFit}) illustrates the superiority of the proposed method in terms of performance. From the single i/o data sequence, a recurrent NN with $3$-layers and $20$-lags was trained using the same training data\footnote{\textbf{Data detrend}:~\footnotesize$u_n=(u-\bar{u})/\sigma_u,~y_n=(y-\bar{y})/\sigma_y$; \textbf{zero-response}: data were doubled in size for learning the zero-response i.e., $u_n=0\rightarrow\boxed{\bSigma}\rightarrow y_n=0$.} as in \cite{Ramos2009} (1000 points). Further, the trained NN was used in the bilinear realization algorithm to generate more data, and a reduced stable bilinear model of order $r=3$ shown in (\ref{eq:daisyfit}) was successfully constructed. The original noisy data were explained with a lower mean percentage error $\text{MPE}=0.56\%$ compared to the subspace method for the whole dataset. Another NN architecture, s.a., the NARMAX\footnote{\textbf{NARMAX}: Nonlinear AutoRegressive Moving Average model with eXogenous inputs \cite{Billings,AlBaiyat2004}.} belongs to a subclass of bilinear systems, and it will filter some nonlinear features without achieving such a good MPE.
\end{example}
In conclusion, NN architectures are a super-class of NARMAX models used in the classical identification theory. Consequently, NN models share the same strong argument with the Carleman linearization scheme that can approximate general nonlinear systems. Finally, NN and realization theory successfully bridge data science with computational science to build reliable, interpretable nonlinear models. Different NN architectures (s.a., recurrent NNs) in combination with other realization frameworks (s.a., Loewner) and for other types of nonlinearities (s.a., quadratic-bilinear) are left for future research endeavors.
\vspace{-5mm}
\addcontentsline{toc}{section}{References}
\bibliographystyle{spmpsci}
\bibliography{main_Karachalios}
\end{document}